\documentclass[12pt]{amsart}
\usepackage{amssymb, latexsym,mathtools}
\usepackage{multicol}
\usepackage{amscd}
\usepackage{extarrows}
\theoremstyle{plain}
\newtheorem{theorem}{Theorem}[section]

\theoremstyle{definition}
\newtheorem{definition}[theorem]{Definition}
\newtheorem{remark}[theorem]{Remark}

\newtheorem{example}[theorem]{Example}

\numberwithin{equation}{section}

\usepackage{dsfont}
\usepackage{braket}
\usepackage{hyperref}

\begin{document}

\title[Barycentric algebras -- convexity and order]{Barycentric algebras - convexity and order}

\author[A. Zamojska-Dzienio]{Anna Zamojska-Dzienio}
\address{Faculty of Mathematics and Information Science\\
Warsaw University of Technology\\
00-661 Warsaw, Poland}

\email{anna.zamojska@pw.edu.pl}


\maketitle

Originating from the use of barycentric coordinates in geometry by M\"obius \cite{Moebius}, barycentric algebras were introduced in the nineteen-forties \cite{Kneser, Stone} for the axiomatization of real convex sets, presented algebraically with binary operations given by weighted means, the weights taken from the (open) unit interval in the real numbers.

Barycentric algebras unify ideas of convexity and order. They provide intrinsic descriptions of convex sets (cancellative b.a.), semilattices (commutative and associative b.a. with all operations being equal), and more general semilattice-ordered systems of convex sets, independently of any ambient affine or vector space. Natural applications can be found in the modeling of systems that function on (potentially incomparable) multiple levels in biology \cite{Sm13} and hierarchical statistical mechanics \cite[Ch.~IX]{Modes}, in theoretical computer science for both nondeterministic and probabilistic
systems for verification \cite{BSV}, in handling thermostatic systems in many contexts \cite{BLM}, and also in computational geometry to analyze systems of barycentric coordinates \cite{RSZ1,RSZ2}.

In this minicourse\footnote{This is the abstract of a series of lectures given during the XIIIth School on Geometry and Physics,
Bia{\l}ystok (Poland), in July 2024.}, we first examine the algebraic aspects of barycentric algebras. Then, we focus on various examples and applications, reviewing the pertinence of the barycentric algebra structure.

\section{Introduction}

In this lecture we will consider (abstract) algebras understood as pairs $(A,F)$ with a set $A$ together with a family $F$ of operations on $A$, i.e. of functions $f\colon A^k\to A$ for some $k\in\mathbb{N}$. In this setting, groups can be understood as algebras with one binary, one unary and one nullary operation (respectively: \emph{multiplication}, taking the \emph{inverse} and the \emph{neutral element}), rings as algebras with two binary, one unary and one nullary operation (rings with unity have two nullary operations), and vector spaces over a field $\mathbb{F}$ as algebras with the group operations and the cardinality of $\mathbb{F}$ many unary operations of \emph{scalar multiplication}. But we focus here on algebras which correspond to affine spaces and convex sets over subfields of $\mathbb{R}$, and barycentric algebras. Our main references are \cite{Modals, Modes, Sm11}.

\section{Affine spaces of $\mathbb{R}^n$}\label{S:as}

Consider the vector space $\mathbb{R}^n$. The line $L_{\mathbf a,\mathbf b}$ through $\mathbf a, \mathbf b \in \mathbb{R}^n$ can be described as
$$
L_{\mathbf a,\mathbf b} = \{\underline{p}(\mathbf a,\mathbf b):= (1-p)\cdot \mathbf a + p\cdot \mathbf b \in \mathbb{R}^n \mid p \in \mathbb{R}\}.
$$
A subset $A \subseteq \mathbb{R}^n$ is a (non-trivial) affine subspace of $\mathbb{R}^n$ if, together with any two distinct points $\mathbf a$ and $\mathbf b$, it contains the line $L_{\mathbf a,\mathbf b}$.

This observation allows us to consider the affine space $A$ as the (abstract) algebra $(A,\underline{\mathbb{R}})$, where $\underline{\mathbb{R}} = \{\underline{p} \mid p \in \mathbb{R}\}$ is a set of continuum many binary operations that are indexed by $\mathbb{R}$.

It is known that affine spaces defined in this way have all the properties required for affine spaces defined in the traditional way. (See e.g. \cite[Sec.~1.6]{Modes}.)
In particular, the affine subpaces of an affine space $A$ are precisely the subalgebras of $(A,\underline{\mathbb{R}})$ (subsets of $A$ closed under the operations $\underline{p}$); the affine transformations $h\colon A \rightarrow A'$ are the algebra homomorphisms $h\colon (A,\underline{\mathbb{R}}) \rightarrow (A',\underline{\mathbb{R}})$ (mappings preserving the operations $\underline{p}$); the direct product of a family of affine spaces is its cartesian product with operations $\underline{p}$ defined componentwise. An $n$-dimensional affine space is affinely isomorphic to the affine space $\mathbb{R}^n$. Elements of the affine hull of a set $\{\mathbf a_1,\dots,\mathbf a_n\}$ of affinely independent elements are described as affine combinations
\begin{equation*}
\mathbf a = \sum_{i=1}^{n} r_i \mathbf a_{i} \ \text{with} \ \sum_{i=1}^{n}r_{i} = 1,
\end{equation*}
where $r_i \in \mathbb{R}$.
The affine combinations may be obtained by composing basic binary operations $\underline{p}$~\cite[Lem.~1.6.4.1, Cor.~6.3.5]{Modes}.
The class $\mathbf{A}(\mathbb{R})$ of (algebras isomorphic to) affine spaces considered as algebras $(A,\underline{\mathbb{R}})$ is a variety. (See e.g. \cite[Sec.~6.3]{Modes}.) This means that $\mathbf{A}(\mathbb{R})$ is closed under the formation of subalgebras, direct products and homomorphic images. As each variety of similar algebras may be equivalently described  as an equationally defined class of algebras, we have the following.

\begin{theorem}\cite{OS}\label{T:OS}
The variety $\mathbf{A}(\mathbb{R})$ is the class of algebras satisfying the identities:
\begin{enumerate}
\item[$(\mathrm a)$] $\underline{p}(x,x)=x$;
\item[$(\mathrm b)$] $\underline{p}\left(\underline{r}(x,y),\underline{r}(z,t)\right)=\underline{r}\left(\underline{p}(x,z),\underline{p}(y,t)\right)$;
\item[$(\mathrm c)$] $\underline{0}(x,y)=x=\underline{1}(y,x)$;
\item[$(\mathrm d)$] $\underline{r}(\underline{p}(x,y),\underline{q}(x,y))=\underline{\underline{r}(p,q)}(x,y)$.
\end{enumerate}
\end{theorem}

In other words, for each algebra $(A,\underline{\mathbb{R}})$ in $\mathbf{A}(\mathbb{R})$, the operations are required to satisfy the properties $(\mathrm a)$-$(\mathrm d)$ of Theorem \ref{T:OS} after substituting any choice of elements of $A$ for the variables $x,y,z,t$.

\section{Convex subsets of $\mathbb{R}^n$}\label{S:cs}

Denote by $I$ the closed real unit interval $[0,1]$, and by $I^{\circ}=I\smallsetminus\{0,1\}=]0,1[$ the open real unit interval. The line segment $I_{\mathbf a,\mathbf b}$ joining the points $\mathbf a,\mathbf b$ of $\mathbb{R}^n$ can be described as
$$
I_{\mathbf a,\mathbf b} = \{\underline{p}(\mathbf a,\mathbf b) \mid p \in I^{\circ}\}.
$$
 A subset $C \subseteq \mathbb{R}^n$ is a (non-trivial) convex subset of $\mathbb{R}^n$ if, together with any two different points $\mathbf a$ and $\mathbf b$, it contains the line segment $I_{\mathbf a,\mathbf b}$. Thus, convex sets can be considered as algebras $(C,\underline{I}^{\circ})$, where $\underline{I}^{\circ}= \{\underline{p} \mid p \in I^{\circ}\}$. They form the class of \emph{subreducts} | subalgebras of the reduct $(A,\underline{I}^{\circ})$ | of affine spaces $(A, \underline{\mathbb{R}})$.

 \begin{remark}\label{R:closed}
 Note that some authors prefer to use the closed unit interval as the indexing set for the basic operations when defining convex sets, instead of the open unit interval. From the standard algebraic point of view, these definitions are equivalent. However, our choice of operations is essential for the algebraic structure theory of barycentric algebras which would not be possible if we also admitted the projections $\underline{0}$ and $\underline{1}$ as basic operations. We will return to this subtle and important point later.
 \end{remark}

 The class $\mathbf{C}$ of (algebras isomorphic to) convex sets considered as algebras $(C,\underline{I}^{\circ})$ is closed under the formation of subalgebras and direct products, but is not closed under homomorphic images (see Example \ref{E:homom}). So, it is not a variety. The class of homomorphic images of convex sets is a variety, the variety $\mathbf{B}$ of \emph{barycentric algebras}.

\section{Barycentric algebras}\label{S:ba}

The class $\mathbf{C}$ and the variety $\mathbf{B}$ can be generalized by taking any subfield of the field $\mathbb{R}$. Thus, we fix a subfield $R$ of $\mathbb{R}$ and take $I^{\circ}=]0,1[\,\cap\, R$ as the open unit interval in $R$. In what follows, we consider algebras with binary operations that are indexed by such $I^\circ$.

\begin{definition}\cite{Modals,Modes,Sm11}\label{D:bar}
The variety $\mathbf{B}$ of barycentric algebras over $R$
is axiomatized by
the identities
\begin{equation*}
\underline{p}(x,x)=x
\end{equation*}
of {\it idempotence},
\begin{equation*}
\underline{p}(x,y)=\underline{1-p}(y,x)
\end{equation*}
of {\it skew-commutativity}, and
\begin{equation*}
\underline{p}(\underline{r}(x,y),z)=\underline{r\circ p}(x,\underline{p/(r\circ p)}(y,z))
\end{equation*}
of {\it skew-associativity}, where $r \circ p = r+p-rp$.
\end{definition}

The variety $\mathbf{B}$, like the variety $\mathbf{A}(\mathbb{R})$, also models the identities of \emph{entropicity} described in Theorem \ref{T:OS}$(\mathrm b)$. One says that the binary operations in $\underline{I}^\circ$ (or in $\underline{\mathbb{R}}$ for affine spaces) \emph{commute}.

The class $\mathbf{C}$ of convex sets (over $R$) is defined, within the variety $\mathbf{B}$, by
the quasi-identities of {\it cancellativity}
\begin{equation*}
\underline{p}(x,y)=\underline{p}(x,z) \mathrel{\Rightarrow} y=z
\end{equation*}
for all $p \in I^{\circ}$.
It has been already shown in \cite[Lemma 2]{N70} (see also \cite[Proposition 5.8.7]{Modes}) that if cancellativity holds in a barycentric algebra for some operation $\underline{p}$, then it holds for all operations in $\underline{I}^{\circ}$.
Note also that all identities defining the variety $\mathbf{B}$ are \emph{regular}, i.e. the sets of variables appearing on both sides of equality symbol $=$ coincide. It would not be the case for operations of projections $\underline{0}$ and $\underline{1}$. (Compare Remark~\ref{R:closed}.)

The idea of abstract convexity has a long history, and has been rediscovered independently many times, using various names for barycentric algebras such as (abstract) convex sets, convex modules, convexors, semiconvex sets, convex spaces, and so on. For historical notes and references, consult \cite[Remarks 2.9]{KP} and \cite[pp. 308-309]{Modes}. More information on real barycentric algebras, i.e. for $R=\mathbb{R}$, is included in \cite{ZD} in the current volume.

\subsection{Examples of barycentric algebras}

We now present some examples crucial for our purposes. An extensive list of other examples can be found in \cite{Fritz} (using the name ``convex space'' for a barycentric algebra) and, for different types of \emph{function spaces}, in \cite{ZD}.

\begin{example}\label{E:vspace}
Let $R=\mathbb{R}$ and $V$ be a vector space over $\mathbb{R}$. Let $\underline{p}\colon V\times V\rightarrow V$ for $p\in I^{\circ}$ be the \emph{weighted mean operation}
$$
\underline{p}(u,v)=(1-p)\cdot u+p\cdot v.
$$
We obtain the algebra $(V,\underline{I}^{\circ})$, and its subalgebras are the convex subsets of the real vector space $V$. For a subfield $R$ of $\mathbb{R}$, the algebra $(V,\underline{I}^{\circ})$ and its subalgebras are cancellative barycentric algebras over $R$.
\end{example}

\begin{example}\label{E:itsemil}
Let $(A,\vee)$ be a (join) semilattice, i.e. $\vee$ is a binary operation satisfying the properties of idempotence, commutativity and associativity. Equivalently, $(A,\vee)$ can be treated as a partially ordered set $(A,\leq)$ in which any two elements $a,b$ have a least upper bound $a\vee b$, so that $a\vee b=b \mathrel{\Leftrightarrow} a\leq b$.
A semilattice $(A,\leq)$ becomes a barycentric algebra $(A,\underline{I}^{\circ})$ (over $R$) if one defines
\begin{equation*}
\underline{p}(a,b)=a\vee b
\end{equation*}
for all $\underline{p}\in \underline{I}^{\circ}$. This is an \emph{iterated semilattice}.
By \cite[Theorem 3]{N70}, iterated semilattices form the only non-trivial proper subvariety $\mathbf{S}$ of $\mathbf{B}$. This is the variety defined by
$$
\underline{p}(x,y)=\underline{r}(x,y)
$$
for all $\underline{p},\underline{r}\in \underline{I}^{\circ}$.
\end{example}

\begin{example}\label{E:homom}
Again, $R=\mathbb{R}$. Consider $(I,\underline{I}^{\circ})\in \mathbf{C}$ (subalgebra of $(\mathbb{R},\underline{I}^{\circ})$) and $(A=\{a,b,c\},\underline{I}^{\circ})\in \mathbf{S}$:
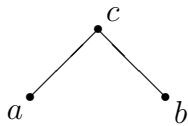
\begin{multicols}{2}
\setlength{\unitlength}{1mm}
\begin{picture}(30,20)(-10,0)
\put(2,2){\makebox(0,0){$a$}} \put(4,4){\circle*{1}}
\put(24,2){\makebox(0,0){$b$}} \put(22,4){\circle*{1}}
\put(15,15){\makebox(0,0){$c$}} \put(13,13){\circle*{1}}
\put(4,4){\line(1,1){9}}
\put(22,4){\line(-1,1){9}}
\end{picture}

  $$\begin{array}{c|ccc}
   \vee & a & b & c\\
   \hline
   a & a & c & c\\
   b & c & b & c\\
   c & c & c & c
   \end{array}$$
\end{multicols}
Take the mapping $h\colon (I,\underline{I}^{\circ})\rightarrow(A,\underline{I}^{\circ})$, $h(0)=a,\; h(1)=b,\; h(x)=c$ for $x\in I^{\circ}$. Direct calculations show that
$h$ is a homomorphism of barycentric algebras
$$
h(\underline{p}(x,y))=\underline{p}(h(x),h(y))=h(x)\vee h(y).
$$
However, $h(I)$ is not a convex set.
\end{example}

Examples \ref{E:vspace} and \ref{E:itsemil} show that in the variety $\mathbf{B}$, defined for some fixed subfield $R$ of $\mathbb{R}$, one can find algebras that behave in distinctly different ways. This is why we say that barycentric algebras unify ideas of convexity and order. Moreover, we can divide them into three disjoint classes: cancellative barycentric algebras (or of \emph{geometric type} as described in \cite{Fritz}), iterated semilattices (or of \emph{combinatorial type} \cite{Fritz}), and the remaining ones (of \emph{mixed type} \cite{Fritz}). Note that in barycentric algebras of mixed type, none of the operations can be \emph{fully} cancellative (by \cite[Lemma 2]{N70}) or a semilattice operation. However, such algebras are constructed by means of convex sets and semilattices, as certain disjoint sums of convex sets. We provide now two important examples of barycentric algebras of mixed type.

\begin{example}\label{E:infty}
The disjoint union of $\mathbb R$ and a singleton $\{\infty\}$ forms a barycentric algebra $(\mathbb R^\infty,\underline{I}^{\circ})$, the \emph{extended real line}, under the operations defined as follows:  $(\mathbb R,\underline{I}^{\circ})$ is a subalgebra of $(\mathbb R^\infty,\underline{I}^{\circ})$, and $\underline p(x,\infty) = \infty$ for all $x\in\mathbb R^\infty$ and $p\in I^\circ$.
\end{example}

\begin{example}\label{E:Talg} Let $T$ denote the union of two copies of $(I,\underline{I}^{\circ})$: $[\alpha,\beta]$ and $[\underline{1/2}(\alpha,\beta),\gamma]$. Suppose that $
[\alpha,\beta]\cap [\underline{1/2}(\alpha,\beta),\gamma]=\{\underline{1/2}(\alpha,\beta)\}$ (see Figure~\ref{F:T}) and define
\begin{equation*}
\underline{p}(x,y)=(1-p)\cdot \underline{1/2}(\alpha,\beta)+p\cdot y
\end{equation*}
for $x\in [\alpha,\beta]$ and $y\notin [\alpha,\beta]$. Then $(T,\underline{I}^{\circ})$ forms a barycentric algebra generating a quasivariety that covers $\mathbf{C}$ in the lattice of subquasivarieties of $\mathbf{B}$ (\cite{Ig} for $\mathbb{R}$, and \cite[Section 7.6]{Modes} for $R$).
\setlength{\unitlength}{1mm}
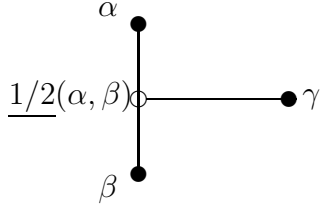
\begin{figure}[hbt]
\begin{picture}(32,26)(-4,0)
\put(3,2){\makebox(0,0){$\beta$}} \put(7,4){\circle*{2}}
\put(-2,14){\makebox(0,0){$\underline{1/2}(\alpha,\beta)$}} \put(7,14){\circle{2}}
\put(3,26){\makebox(0,0){$\alpha$}} \put(7,24){\circle*{2}}
\put(30,14){\makebox(0,0){$\gamma$}} \put(27,14){\circle*{2}}
\put(7,4){\line(0,1){20}}
\put(8,14){\line(1,0){20}}
\end{picture}
\caption{The barycentric algebra $(T,\underline{I}^{\circ})$}
\label{F:T}
\end{figure}
\end{example}

\section{The structure of barycentric algebras}\label{S:ss}

The main purpose of this part is to present some theorems describing the structure of general
barycentric algebras. They formulate precisely the way barycentric algebras combine order and
convexity. This section is based mainly on \cite{RS90} and \cite[Section 7.5]{Modes}.

\subsection{Semilattice sums}

By \cite[Theorem 3]{N70} we know that each bary\-centric algebra has a homomorphism onto an iterated semilattice. The largest semilattice quotient (homomorphic image) $S$ of a barycentric algebra $A$ is called its \emph{semilattice replica}. Let $\varrho:A \rightarrow S$ be the surjective \emph{semilattice replica homomorphism} of $A$ witnessing the semilattice replica $S$. For each $s \in S$, the preimage $\varrho^{-1}(\{s\})$ is a subalgebra of $A$ and a convex set. In such a context, we will say that the barycentric algebra is a \emph{semilattice sum of convex sets}. (For details see \cite[Section 7.5]{Modes}.) We obtain the following.

\begin{theorem}\label{T:strthm}\cite{RS90}\cite[Corollary 7.5.9]{Modes}
Each barycentric algebra is a semilattice sum of open convex sets.
\end{theorem}

\begin{remark}\label{R:wall}
\begin{enumerate}
\item[$(\mathrm a)$] For a barycentric algebra $(A,\underline{I}^{\circ})$, a subset $W$ of $A$ is a \emph{wall}, if
$$
\forall a,b\in A\;\forall r \in I^{\circ}
\;\;\underline{r}(a,b) \in W\mathrel{\Leftrightarrow} a\in W \text{ and } b\in W.
$$
Note that in particular, according to the implication $\Leftarrow$, walls are special kinds of subalgebras of $(A,\underline{I}^{\circ})$. An (\emph{algebraically}) \emph{open} barycentric algebra has no proper non-empty walls.
\item[$(\mathrm b)$] The semilattice underlying Theorem \ref{T:strthm} is the semilattice replica of the barycentric algebra. The crucial property here and in other structure theorems based on Theorem \ref{T:strthm} is that the identities satisfied by barycentric algebras are regular. (Recall Remark \ref{R:closed}.)
\end{enumerate}
\end{remark}

\begin{example}\label{E:semsum}
\begin{enumerate}
\item[$(\mathrm a)$] Return to Example \ref{E:homom}. The iterated semilattice $(A,\underline{I}^{\circ})$ is the semilattice replica of $(I,\underline{I}^{\circ})$ with open convex sets $\varrho^{-1}(\{a\})=\{0\}$, $\varrho^{-1}(\{b\})=\{1\}$ and $\varrho^{-1}(\{c\})=I^{\circ}$.
\item[$(\mathrm b)$] The semilattice replica of the extended real line (Example \ref{E:infty}) is the two-element semilattice $(\{a,b\},\underline{I}^{\circ})$, with $a<b$, and open convex sets $\varrho^{-1}(\{a\})=\mathbb{R}$, $\varrho^{-1}(\{b\})=\{\infty\}$.
\item[$(\mathrm c)$] The semilattice replica of the algebra $(T,\underline{I}^{\circ})$ (Example \ref{E:Talg}) is the five-element semilattice $(\{a,b,c,d,e\},\underline{I}^{\circ})$ (see Figure \ref{F:T_repl}). Now, in its decomposition, we obtain the open convex sets $\varrho^{-1}(\{a\})=\{\alpha\}$, $\varrho^{-1}(\{b\})=\{\beta\}$, $\varrho^{-1}(\{c\})=]\alpha,\beta[$, $\varrho^{-1}(\{d\})=\{\gamma\}$ and $\varrho^{-1}(\{e\})=]\underline{1/2}(\alpha,\beta),\gamma[$.
\end{enumerate}
\end{example}

\setlength{\unitlength}{1mm}
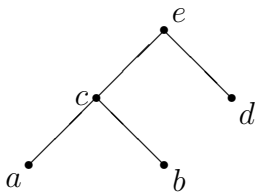
\begin{figure}[hbt]
\begin{picture}(32,26)(0,0)
\put(2,2){\makebox(0,0){$a$}} \put(4,4){\circle*{1}}
\put(24,2){\makebox(0,0){$b$}} \put(22,4){\circle*{1}}
\put(11,13){\makebox(0,0){$c$}} \put(13,13){\circle*{1}}
\put(33,11){\makebox(0,0){$d$}} \put(31,13){\circle*{1}}
\put(24,24){\makebox(0,0){$e$}} \put(22,22){\circle*{1}}
\put(4,4){\line(1,1){9}}
\put(22,4){\line(-1,1){9}}
\put(13,13){\line(1,1){9}}
\put(22,22){\line(1,-1){9}}
\end{picture}
\caption{The semilattice replica of $(T,\underline{I}^{\circ})$}
\label{F:T_repl}
\end{figure}

\subsection{P{\l}onka sums}

Theorem~\ref{T:strthm} provides a decomposition of a bary\-centric algebra $(A,\underline{I}^{\circ})$ into disjoint open convex subsets $(A_s,\underline{I}^{\circ})$ over its semilattice replica $S$. But we would also like to have a construction method that allows to recover the algebra $(A,\underline{I}^{\circ})$ from the fibres $A_s$, certain extensions $E_s$ of $A_s$, and the semilattice $S$. This construction was described in \cite{RS90}. The extensions $E_s$ are also convex sets, in a sense minimal, and are obtained  as homomorphic images of some special subalgebras of $(A,\underline{I}^{\circ})$. There is a canonical way to extend the structure of a barycentric algebra on $(A,\underline{I}^{\circ})$ to the disjoint sum $E$ of $E_s$, where $s \in S$, in such a way that $E$ is a semilattice sum of $E_s$ over $S$ with nice properties. This type of semilattice sum was introduced by P\l onka in \cite{Pl} as a generalization of the \emph{strong semilattice of semigroups} to arbitrary algebras. For a contemporary approach to \emph{P{\l}onka sums} see \cite[Chapter 4]{Modes} or \cite[Chapter 2]{Bonzio}.

First recall that a (join) semilattice $(S,\vee)$ may be considered as a category, with the elements of $S$ as objects, and with morphisms $s \rightarrow t$ if $s \leq t$.
We may also consider the quasivariety $\mathbf{C}$ as a category, with convex sets as objects and homomorphisms as morphisms. Let $F: S \rightarrow \mathbf{C}$ be a (covariant) functor from the category $S$ to the category $\mathbf{C}$ which acts on morphisms as follows:
$$
\begin{array}{cccc}
&t\,\,  &\phantom{-}  &(C_t,\underline{I}^{\circ})\\
&\Big\uparrow  &\xlongrightarrow{\;F\; } &\Big\uparrow \varphi_{s,t}\\
&s\,\,  &\phantom{-}  &(C_s,\underline{I}^{\circ})
\end{array}
$$

\begin{definition}
Let $F:S \rightarrow \mathbf{C}$ be a functor from a semilattice $(S,\vee)$ to the category $\mathbf{C}$ of convex sets.
The \emph{P\l onka sum} of convex sets $(C_s,\underline{I}^{\circ})$ over the semilattice $(S,\vee)$ by the functor $F$ is the algebra $(A, \underline{I}^{\circ})\in\mathbf{B}$ defined on the disjoint sum $\biguplus_{s\in S}C_{s}$ of $C_s$, with the operations $\underline{p} \in \underline{I}^{\circ}$ given by
$$
\underline{p}\colon C_{s} \times C_{t} \rightarrow C_{s\vee t}; (a_{s}, a_{t}) \mapsto \underline{p}(\varphi_{s,s\vee t}(a_s )\,, \varphi_{t, s\vee t}(a_t)).
$$
\end{definition}

The meaning of the P\l onka sum construction for barycentric algebras is given by the following Structure Theorem.

\begin{theorem}\cite{RS90}\label{T:str}
Each barycentric algebra is a subalgebra of a P\l onka sum
of convex sets over its semilattice replica.
\end{theorem}

\begin{example}\label{E:Plsum}
\begin{enumerate}
\item[$(\mathrm a)$] The interval $(I,\underline{I}^{\circ})$ from Example \ref{E:homom} is a subalgebra of the P\l onka sum of convex sets $E_a=\{0\}$, $E_b=\{1\}$ and $E_c=I$ over the iterated semilattice $(\{a,b,c\},\underline{I}^{\circ})$. The P\l onka homomorphism $\varphi_{a,c}$ maps $0$ to the left end of $I$ and, similarly, $\varphi_{b,c}$ maps $1$ to the right end.
\item[$(\mathrm b)$] The semilattice sum of the barycentric algebra $(\mathbb R^\infty,\underline{I}^{\circ})$ considered in Examples \ref{E:infty} and \ref{E:semsum}$(\mathrm b)$ is already a P\l onka sum of  $E_a=\mathbb{R}$, $E_b=\{\infty\}$.
\item[$(\mathrm c)$] The algebra $(T,\underline{I}^{\circ})$ (Example \ref{E:Talg}) is the P\l onka sum of convex sets $E_a=\{\alpha\}$, $E_b=\{\beta\}$, $E_c = [\alpha,\beta]$, $E_d=\{\gamma\}$ and $E_e = [\underline{1/2}(\alpha,\beta),\gamma]$. The P\l onka homomorphism $\varphi_{c,e}$ maps the whole of $E_c$ into the element $\underline{1/2}(\alpha,\beta)$ of $E_e$. The remaining P\l onka homomorphisms behave as in the case $(\mathrm a)$ above.
\end{enumerate}
\end{example}

\subsection{Semilattice sums for convex polytopes}

Consider now \emph{convex polytopes}. As algebras, they are defined as finitely generated convex sets. The minimal set of generators of a polytope is its set of vertices (i.e., its extreme points). In geometric terminology, the convex set generated by a set $V$ is its \emph{convex hull}. If a $k$-dimensional polytope has $n$ vertices, then $n$ is at least $k+1$.

One can treat a convex polytope $\mathcal{C}$ in $\mathbb{R}^k$ as a cancellative barycentric algebra $(C,\underline{I}^{\circ})$, a subalgebra of $(\mathbb{R}^k,\underline{I}^{\circ})\in\mathbf{C}$ constructed as described in Example~\ref{E:vspace}. The set of vertices $V$ of $\mathcal{C}$ is the generating set of the algebra $(C,\underline{I}^{\circ})$, i.e. $(C,\underline{I}^{\circ})$ is the smallest cancellative barycentric algebra which contains $V$.

Note that in the case of convex polytopes, the geometric concept of a face and the algebraic concept of a wall (see Remark \ref{R:wall}) coincide. The $0$-dimensional faces of a $k$-dimensional polytope $\mathcal{C}$ are its vertices, and the $1$-dimensional faces are the edges. The maximal (thus $(k-1)$-dimensional) faces are also called the \emph{facets} of $\mathcal{C}$. The faces of a polytope are again polytopes, and each face of $\mathcal{C}$ considered as a barycentric algebra is generated by the vertices of $\mathcal{C}$ that it contains. The faces of a polytope $\mathcal{C}$ form a lattice under inclusion. A polytope $\mathcal{C}$, as a barycentric algebra $(C,\underline{I}^{\circ})$, is a semilattice sum of open convex subsets $C_s$ (trivial faces and interiors of nontrivial faces) over the semilattice replica $S$ that is the join (semi)lattice of its non-empty faces. Moreover, the P\l onka sum $E$ extending a polytope $\mathcal{C}$ has a very simple form. The summands $E_s$ are just the faces of $\mathcal{C}$, and each $C_s$ is the interior of $E_s$. The P\l onka homomorphisms $\varphi_{s,t}: E_s \rightarrow E_t$ for $s \leq t$ are just embeddings (injective homomorphisms). (See Example \ref{E:Plsum}$(\mathrm a)$.)

\subsection{A toy model for complex systems}

	The significance of bary\-centric algebras is that they provide a general algebraic framework for the study of both convexity and order. We have already seen that this happens on many different levels. For example, barycentric algebras are ordered systems of convex sets with order given by a semilattice. Now we will show how this can be applied to the modeling of systems that function on (potentially incomparable) multiple levels. We use a toy model for complex systems, based on the algebra $(T,\underline{I}^{\circ})$ introduced in \cite[Section 16]{Sm11} and described in more detail in \cite[Section 4.8]{Sm13}.

\begin{remark}\label{R:T}
Note that the algebra $(T,\underline{I}^{\circ})$ can be also presented as a subalgebra of the P\l onka sum of two closed intervals $I_0 = [\alpha,\beta]$, and $I_1 = [\underline{1/2}(\alpha,\beta),\gamma]$ over the two-element semilattice with $0<1$. Then the P\l onka homomorphism $\varphi_{0,1}$ maps the whole of $I_0$ into $\underline{1/2}(\alpha,\beta)$ of $I_1$. This is due to the fact that there exists a semilattice homomorphism from the semilattice replica to the two-element semilattice given by $a,b,c\mapsto 0$ and $d,e\mapsto 1$.
\end{remark}

We consider a model in biology. Suppose there are two species, $A$ and $B$. Species $A$ exists in two stages, $A1$ and $A2$ (e.g. larva and adult), while species $B$ is unstructured. Now, $A$ and $B$ compete for a limited supply of food, such that the total number of individuals in $B$ and $A$ (regardless of stage) is constant. We have two levels here: \emph{demography}, which deals with the internal structure of a single species, and \emph{ecology} which addresses the competition between different species. The system is modeled by the algebra $(T,\underline{I}^{\circ})$ (see Figure \ref{F:toy}). Demographic states (mixes of $A1$ and $A2$) are elements of the subalgebra $([A1,A2],\underline{I}^{\circ})$, while ecological states correspond to elements of  $([\underline{1/2}(A1,A2),B],\underline{I}^{\circ})$. The biology here says it does not matter in which stage a particular individual of $A$ exists, so at the ecological level, each particular individual of $A$ is represented by a uniform mix $\underline{1/2}(A1,A2)$.
\setlength{\unitlength}{1mm}
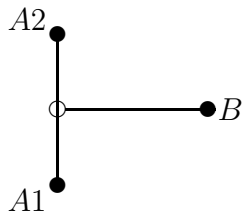
\begin{figure}[hbt]
\begin{picture}(32,26)(-4,0)
\put(3,2){\makebox(0,0){$A1$}} \put(7,4){\circle*{2}}
\put(7,14){\circle{2}}
\put(3,26){\makebox(0,0){$A2$}} \put(7,24){\circle*{2}}
\put(30,14){\makebox(0,0){$B$}} \put(27,14){\circle*{2}}
\put(7,4){\line(0,1){20}}
\put(8,14){\line(1,0){20}}
\end{picture}
\caption{Demographic and ecological levels}
\label{F:toy}
\end{figure}

\section{The transition from affine to projective geometry}\label{S:trans}
We now want to show how barycentric algebras give a framework for studying affine and projective geometry and the relationship between them. This part is based on \cite{RS85}.

First, we return to affine spaces, but over an arbitrary field $\mathbb{K}$, not necessarily $\mathbb{R}$. We still treat them as (abstract) algebras, now of the form $(V,\underline{\mathbb{K}},P)$ with $V\subseteq\mathbb{K}^n$, a set of binary operations $\underline{\mathbb{K}}=\{\underline{k}\colon k\in \mathbb{K}\}$, where $\underline{k}(u,v)=(1-k)\cdot u+k\cdot v$, and one additional ternary operation $$P\colon V\times V\times V\rightarrow V;(u,v,w)\mapsto u-v+w$$ called the \emph{parallelogram operation}. If the element $2=1+1$ is invertible in $\mathbb{K}$, then $P(u,v,w)=\underline{2}(v,\underline{2^{-1}}(u,w))$, and one can consider affine spaces as algebras $(V,\underline{\mathbb{K}})$ (see \cite{OS} or \cite[Section 1.6.4]{Modes}). In particular, this is the case for $\mathbb{K}=\mathbb{R}$. The class $\mathbf{A}(\mathbb{K})$ of (algebras isomorphic to) affine spaces considered as algebras $(V,\underline{\mathbb{K}},P)$ is a variety whose axiomatization can be found e.g. in \cite{OS}. (Theorem \ref{T:OS} is the version for $\mathbb{K}$ with $2$ invertible.)

The algebra $(V,\underline{\mathbb{K}},P)$ has the affine group as
its group of automorphisms, and may thus be identified as the affine geometry. The corresponding projective geometry consists of the set $\mathcal{L}(V)$ of linear subspaces of $V$ ordered by inclusion, and may be described algebraically as a semilattice $(\mathcal{L}(V),\vee)$ with
$$
U_1\vee U_2:=\{u_1+u_2\colon u_1\in U_1,\; u_2\in U_2\}\in \mathcal{L}(V),
$$
for $U_1,\; U_2\in \mathcal{L}(V)$.
Now take the set $\mathcal{S}(V,\underline{\mathbb{K}},P)$ of (non-empty) subalgebras of $(V,\underline{\mathbb{K}},P)$. One can again introduce the structure of a $\underline{\mathbb{K}}\cup\{P\}$-algebra on $\mathcal{S}(V,\underline{\mathbb{K}},P)$ as follows:
\begin{enumerate}
\item[$(\mathrm i)$] for $k\in\mathbb{K}$: $\underline{k}(U_1,U_2)=\{\underline{k}(u_1,u_2)\colon u_1\in U_1,\; u_2\in U_2\}$
\item[$(\mathrm ii)$] for $P$: $P(U_1,U_2,U_3)=\{P(u_1,u_2,u_3)=u_1-u_2+u_3\colon u_i\in U_i,\; i=1,2,3\}$
\end{enumerate}
due to the idempotence and entropicity of $(V,\underline{\mathbb{K}},P)$.

Now introduce $\Omega_\mathbb{K}\subsetneq \underline{\mathbb{K}}\cup\{P\}$, understood as the algebraic analogue of $I^{\circ}$ in the field of rationals $\mathbb{Q}$. The algebra $(\mathcal{S}(V,\underline{\mathbb{K}},P),\underline{\Omega_\mathbb{K}})$ is a P\l onka sum over the semilattice $(\mathcal{L}(V),\underline{\Omega_\mathbb{K}})$, which follows from the following general proposition.

\begin{theorem}\cite[Proposition 3.4]{RS85}
Let $\mathbb{K}$ be a field other than GF$(2)$ and let $k\in\mathbb{K}\smallsetminus\{0,1\}$. Let $V$ be a vector space over $\mathbb{K}$, with $(\mathcal{L}(V),\vee)$ as the corresponding projective space.
\begin{enumerate}
\item[$(\mathrm i)$] There is a homomorphism:
$$
\pi\colon (\mathcal{S}(V,\underline{\mathbb{K}},P),\underline{k})\rightarrow (\mathcal{L}(V),\vee);\; x+U\mapsto U.
$$
\item[$(\mathrm ii)$] For the functor $F\colon \mathcal{L}(V)\rightarrow(\{\underline{k}\})$ with $F(U)=\pi^{-1}(U)$ and
$$
\begin{array}{cccccc}
&U_t\,\,  &\phantom{-}  &(\pi^{-1}(U_t),\underline{k}) &\phantom{-} &x+U_t\\
&\Big\uparrow  &\xlongrightarrow{\;F\; } &\;\;\Big\uparrow \varphi_{s,t} &\phantom{-} &\rotatebox{90}{$\mapsto$}\\
&U_s\,\,  &\phantom{-}  &(\pi^{-1}(U_s),\underline{k}) &\phantom{-} &x+U_s
\end{array}
$$
the algebra $(\mathcal{S}(V,\underline{\mathbb{K}},P),\underline{k})$ is the P\l onka sum over $(\mathcal{L}(V),\vee)$ by the functor $F$.
\end{enumerate}
\end{theorem}

The main result of \cite{RS85} that gives the direct
invariant passage from affine to projective geometry follows.

\begin{theorem}\cite[Theorem 2.4]{RS85}
Let $(V,\underline{\mathbb{K}},P)$ be an affine space over a field $\mathbb{K}$. Then the projective space $(\mathcal{L}(V),\vee)$ is the semilattice replica of the algebra $(\mathcal{S}(V,\underline{\mathbb{K}},P),\underline{\Omega_\mathbb{K}})$.
\end{theorem}

The proof consists of three cases, depending on the characteristic of the field: $\mathbb{K}$ of characteristic $0$, of odd characteristic, and of characteristic $2$. We present here the sketch of the proof for $\mathbb{K}$ of characteristic $0$ as a good application of what we have learned so far.

We will work with the rational number field $\mathbb{Q}$, as a prime field of $\mathbb{K}$ of characteristic $0$. Here, the element $2$ is invertible. We consider an affine space $(V,\underline{\mathbb{Q}})$ and its reduct $(V,\underline{I}^{\circ})$, where $I^{\circ}=\{x\in\mathbb{Q}\colon 0<x<1\}$. (Note that $\Omega_\mathbb{K}=I^{\circ}$ for $\mathbb{K}$ of characteristic $0$.) We obtain \emph{rational barycentric algebras} (one adjusts Definition \ref{D:bar} to $R=\mathbb{Q}$). Subalgebras of $(V,\underline{I}^{\circ})$ are $\mathbb{Q}$-convex subsets of a vector space $V$.

Now we consider the algebra $(\mathcal{S}(V,\underline{\mathbb{K}}),\underline{I}^{\circ})$ being the $\underline{I}^{\circ}$-reduct of the algebra of affine subspaces of $(V,\underline{\mathbb{K}})$. Its structure is described as follows.

\begin{theorem}\cite[Theorem 3.5]{RS85}\label{T:T1}
The algebra $(\mathcal{S}(V,\underline{\mathbb{K}}),\underline{I}^{\circ})$ is a P\l onka sum of $\mathbb{Q}$-convex sets over the projective space $(\mathcal{L}(V),\vee)$ by the functor $F\colon \mathcal{L}(V)\rightarrow \mathbf{C}$ with $F(U)=(\{x+U\colon x\in V\},\underline{I}^{\circ})$ and
$$
\begin{array}{cccccc}
&U_t\,\,  &\phantom{-}  &F(U_t) &\phantom{-} &x+U_t\\
&\Big\uparrow  &\overset{F}{\longrightarrow} &\Big\uparrow \varphi_{s,t} &\phantom{-} &\rotatebox{90}{$\mapsto$}\\
&U_s\,\,  &\phantom{-}  &F(U_s) &\phantom{-} &x+U_s
\end{array}
$$
\end{theorem}
In particular, we obtain that $(\mathcal{S}(V,\underline{\mathbb{K}}),\underline{I}^{\circ})$ is a rational barycentric algebra (but, again, the crucial property here is that we only have regular identities in our definition of barycentric algebras, see Remark \ref{R:closed}).

Summarizing, we have exhibited the rational barycentric algebra of affine subspaces of an affine $\mathbb{K}$-space as a P\l onka sum of $\mathbb{Q}$-convex sets over a projective space.
In \cite[Theorem 3.8]{RS85} it was also shown that for an affine space $(V,\underline{\mathbb{K}})$, the reduct $(V,\underline{I}^{\circ})$ has no non-trivial $\underline{I}^{\circ}$-semilattice quotients. Hence, the projective geometry $(\mathcal{L}(V),\vee)$ is the largest semilattice quotient of $(\mathcal{S}(V,\underline{\mathbb{K}}),\underline{I}^{\circ})$, i.e. its semilattice replica. (Compare to the situation described in Remark \ref{R:T}.)

\section*{Acknowledgment}
I would like to express my gratitude to Prof. Anna B. Romanowska and
Prof. Jonathan D.H. Smith for introducing me to the wonderful world of
\emph{modes} (idempotent and entropic algebras), and in particular to
barycentric algebras. A nice coincidence is that 2025 is the
40th anniversary of the publication of \cite{RS85} and \cite{Modals}, where the
general investigation of these algebras really got under way.

\end{document}